\documentclass[a4paper]{amsart}
\numberwithin{equation}{section}
\newtheorem{theorem}{Theorem}[section]
\newtheorem{corollary}{Corollary}[section]
\newtheorem{definition}{Definition}[section]
\newtheorem{lemma}{Lemma}[section]

\theoremstyle{remark}

\usepackage{amsmath, hyperref}
\usepackage{color}
\usepackage{graphicx}
\usepackage{subfigure}

\title[A note on the majorization results for certain starlike functions]{A note on the paper "Tang et al. [Bull Iran Math Soc (2019) doi:10.1007/s41980-019-00262-y]"}

\subjclass[2010]{30C45, 30C80}

\keywords{univalent, starlike, majorization, subordination, cosine function, Booth lemniscate.\\
*Corresponding Author, ORCID iD: https://orcid.org/0000-0003-1029-5386\\
Mobile:+989141667438
}

\begin{document}
\begin{abstract}
Very recently Tang et al. [Bull Iran Math Soc (2019) doi:10.1007/s41980-019-00262-y] have studied some majorization results for two certain subclasses of the starlike functions associated with the sine and cosine functions defined by $\mathcal{S}^*_s$ and $\mathcal{S}^*_c$, respectively. In this note we pointed out that the definition of the class $\mathcal{S}^*_c$ and it's result are incorrect and give correct definition and result.
\end{abstract}

\author[Mohammad Mahdi Motamedinezhad and Rahim Kargar] {Mohammad Mahdi Motamedinezhad and Rahim Kargar$^*$}
\address{Department of Mathematics, University of Applied Science and Technology,
Tehran, Iran} \email {motamedi@uast.ac.ir}
\address{Young Researchers and Elite Club,
Ardabil Branch, Islamic Azad University, Ardabil, Iran}
       \email{rkargar1983@gmail.com}

\maketitle

\section{Introduction}
We denote by $\mathcal{H}$ the family of analytic functions $f$ on the open unit disc $\Delta=\{z\in \mathbb{C}:|z|<1\}$ and by $\mathcal{A}\subset\mathcal{H}$ the class of normalized functions $f$ of the form
\begin{equation}\label{f}
f(z)=z+ \sum_{n=2}^{\infty}a_{n}z^{n}\quad(z\in\Delta).
\end{equation}
The subclass of $\mathcal{A}$ consisting of all univalent functions $f(z)$ in $\Delta$ will be denoted by $\mathcal{U}$.
Also we denote by $\mathcal{B}\subset\mathcal{H}$ the family of functions $f$ on $\Delta$ satisfying the conditions $\phi(0)=0$ and $|\phi(z)|<1$. A function $\phi$ belongs to the class $\mathcal{B}$ is said to be a Schwarz function on $\Delta$.

Let two functions $f$ and $g$ belong to the class $\mathcal{H}$ and there exists a Schwarz function $\phi$ such that $f(z)=g(\phi(z))$. Then we say that $f$ is subordinate to $g$, written as $f(z)\prec g(z)$ or $f\prec g$. It is clear that if $f\prec g$, then we have
\begin{equation}\label{conditions subordination}
f(0)=g(0)\quad {\rm and}\quad f(\Delta)\subset g(\Delta).
\end{equation}
In particular, if $g$ belong to the class $\mathcal{U}$, then $f\prec g$ if and only if the conditions \eqref{conditions subordination} hold true. In the sequel, we recall the following definition from \cite{McC}.
\begin{definition}\label{maj def}
  Let $f$ and $g$ belong to the class $\mathcal{H}$. A function $f$ is said to be majorized by $g$ written as $f(z)\ll g(z)$ or $f\ll g$, if there exists an analytic function $\psi$ in $\Delta$, satisfying
  \begin{equation}\label{conditions maj}
    |\psi(z)|\leq 1\quad{\rm and}\quad f(z)=\psi(z)g(z)
  \end{equation}
  for all $z\in\Delta$.
\end{definition}

For an analytic univalent function $\varphi$ with ${\rm Re}\{\varphi(z)\}>0$ $(z\in\Delta)$ and normalized by $\varphi(0)=1$ and $\varphi'(0)>0$, Ma and Minda \cite{Ma1992} introduced the class $\mathcal{S}^*(\varphi)$. A function $f\in\mathcal{A}$ belongs to the class $\mathcal{S}^*(\varphi)$ if and only if $zf'(z)/f(z)$ is subordinate to $\varphi$.
By choosing some special function $\varphi$, several authors have defined many new subclasses of the starlike functions in recent years, see for instance \cite{choBIran, KarEbaSok, KargarComplex, KO2011, Robertson, SokSta}.

Very recently, Tang et al. \cite{Tang2019} introduced the class $\mathcal{S}^*_c$ as follows:
\begin{equation*}
  \mathcal{S}^*_c:=\left\{f\in\mathcal{A}:\frac{zf'(z)}{f(z)}\prec 1+\cos z\right\}.
\end{equation*}
Also, they proved the following theorem:\\
{\bf Theorem A.} {\it Let the function $f\in\mathcal{A}$ and suppose that $g\in\mathcal{S}^*_c$. If $f(z)\ll g(z)$, then for $|z|\leq r_0$
\begin{equation*}
  |f'(z)|\leq|g'(z)|,
\end{equation*}
where $r_0$ is the smallest positive root of the following equation:
\begin{equation*}
  (1-r^2)(1-\cosh r)-2r=0.
\end{equation*}}
First of all we notice that the definition of the class $\mathcal{S}^*_c$ is incorrect, because if $f\in\mathcal{A}$, then
\begin{equation*}
  \frac{zf'(z)}{f(z)}\Big|_{z=0}=1
\end{equation*}
and $(1+\cos z)|_{z=0}=2$ while $(zf'(z)/f(z))\prec (1+\cos z)$.
Therefore Theorem A is incorrect, too. In this note, we first will define the class $\mathcal{S}^*_c$ correctly and obtain majorization property (the correct version of Theorem A) for it. 

In order to prove the main result, we need the following lemma, see \cite{Nehari1952}.
\begin{lemma}\label{Nehari lemma}
  Let $\psi(z)$ be analytic in $\Delta$ and satisfying $|\psi(z)|\leq1$ for all $z\in\Delta$. Then
  \begin{equation}\label{estimate Nehari}
    |\psi'(z)|\leq\frac{1-|\psi(z)|^2}{1-|z|^2}.
  \end{equation}
\end{lemma}
\section{Main Result}
First, we recall that if $f(z)$ is majorized by $g(z)$ in $\Delta$ and $g(0)=0$, then
\begin{equation*}
  \max_{|z|=r}|f'(z)|\leq \max_{|z|=r}|g'(z)|
\end{equation*}
for each number $r$ in the interval $[0,\sqrt{2}-1]$, see \cite{McC}. By choosing $g(z)=z$, the above assertion generalizes the following well-known theorem \cite{Cara}:\\
{\bf Theorem B.} {\it If the function $f\in\mathcal{H}$ satisfies the conditions
\begin{equation*}
  |f(z)|\leq 1\quad{and}\quad f(0)=0,
\end{equation*}
then $|f'(z)|\leq 1$ for $|z|\leq \sqrt{2}-1$.}

\begin{figure}[!ht]
\centering
\subfigure[]{
\includegraphics[width=4.5cm]{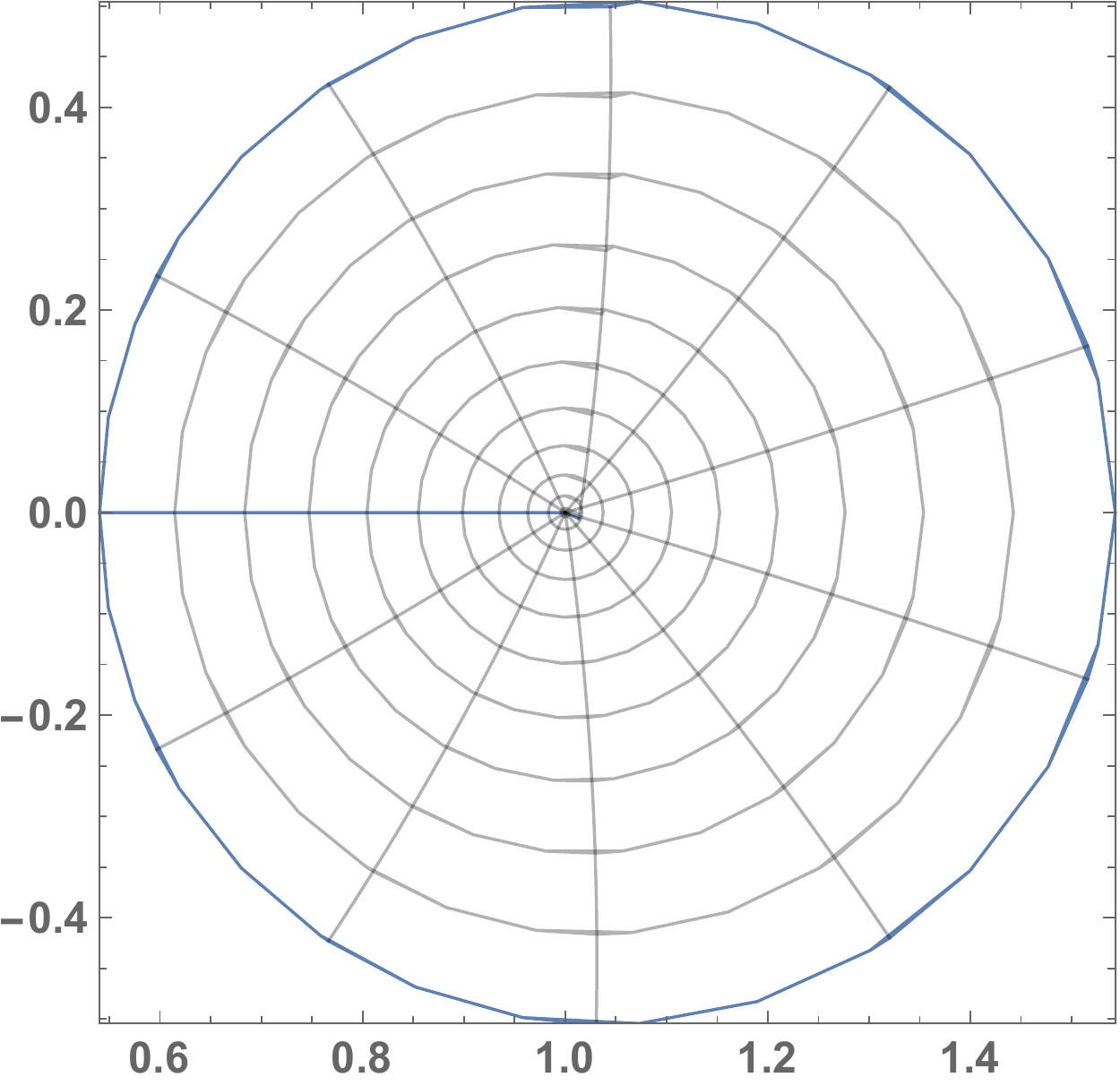}
	\label{fig:subfig1}
}
\hspace*{10mm}
\subfigure[]{
\includegraphics[width=5.5cm]{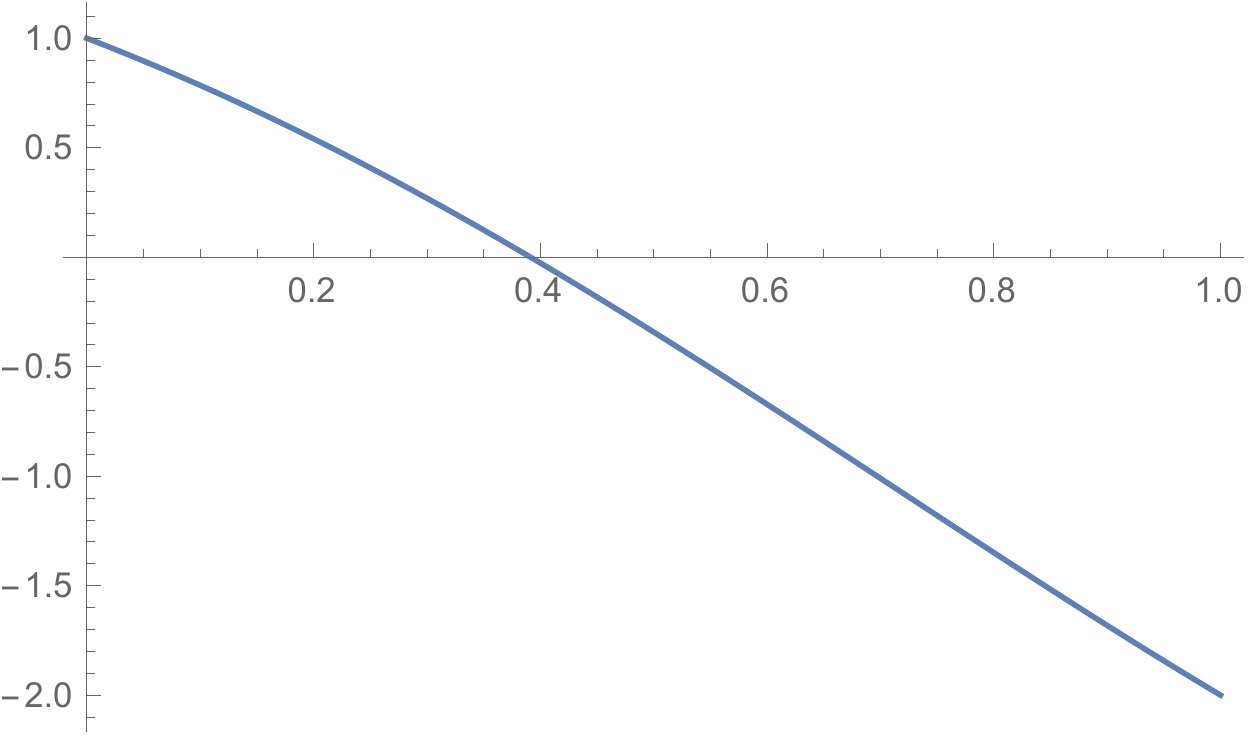}
	\label{fig:subfig2}
}

\caption[The boundary curve of $\cos(\Delta)$]
{\subref{fig:subfig1}: The boundary curve of $\cos(\Delta)$  ،
 \subref{fig:subfig2}: The graph of $(1-r^2)\cos r-2r=0$ when $0<r<1$،
 }
\label{fig:subfig01}
\label{Fig1}
\end{figure}

Following we define the correct version of the class $\mathcal{S}^*_c$ and denote by $\mathbf{S}^*_c$.
\begin{definition}\label{def of the correct version of S c}
  We say that a function $f\in\mathcal{A}$ belongs to the class $\mathbf{S}^*_c$ if and only if
\begin{equation}
  \frac{zf'(z)}{f(z)}\prec \cos z.
\end{equation}
\end{definition}
Since $\cos z$ is univalent in $\Delta$ and $\cos z|_{z=0}=1$, thus the above Definition \ref{def of the correct version of S c} is well-defined. Figure \ref{fig:subfig1} shows the image of $\Delta$ under the function $\cos z$.

\begin{theorem}\label{th1}
  Let $f\in\mathcal{A}$ and $g\in\mathbf{S}^*_c$. If $f(z)$ is majorized by $g(z)$ in $\Delta$, then 
  \begin{equation*}
    |f'(z)|\leq |g'(z)|\quad(|z|\leq r_1),
  \end{equation*}
where $r_1\approx 0.391389$ is the smallest positive root of the following equation:
  \begin{equation}\label{eq 1}
    (1-r^2)\cos r-2r=0\quad(0<r<1).
  \end{equation}
\end{theorem}
\begin{proof}
  Let $f\in\mathcal{A}$ be given by \eqref{f} and $g\in\mathbf{S}^*_c$. Thus by Definition \ref{def of the correct version of S c} and by the subordination principle there exists a Schwarz function $\phi$ such that
  \begin{equation*}
    \frac{zg'(z)}{g(z)}= \cos (\phi(z))\quad(z\in\Delta),
  \end{equation*}
or equivalently
\begin{equation}\label{g bar g prime}
  \frac{g(z)}{g'(z)}=\frac{z}{\cos(\phi(z))}\quad(z\in\Delta).
\end{equation}
Let $\phi(z)=Re^{it}$ where $R\leq r=|z|<1$ and $-\pi\leq t\leq \pi$. It is a simple exercise that
\begin{equation}\label{estimate abs cos phi}
  \cos r\leq \cos R\leq |\cos(\phi(z))|\leq \cosh R\leq \cosh r\quad(R\leq r<1).
\end{equation}
From \eqref{g bar g prime} and \eqref{estimate abs cos phi} we get
\begin{equation}\label{estimate 2}
  \left|\frac{g(z)}{g'(z)}\right|\leq \frac{r}{\cos r}\quad(0<r<1).
\end{equation}
Also, since $f$ is majorized by $g$ in $\Delta$, thus we can find an analytic function $\psi$ in $\Delta$ with $|\psi(z)|\leq1$ such that
\begin{equation}\label{f = psi g}
  f(z)=\psi(z)g(z)\quad(z\in\Delta).
\end{equation}
Taking a simple derivation of relation \eqref{f = psi g} we get
\begin{equation}\label{f prime =}
  f'(z)=\psi'(z)g(z)+\psi(z)g'(z)=g'(z)\left(\psi'(z)\frac{g(z)}{g'(z)}+\psi(z)\right).
\end{equation}
If we apply \eqref{estimate 2} and Lemma \ref{Nehari lemma} in \eqref{f prime =}, then we obtain
\begin{equation}\label{abs f prime}
  |f'(z)|\leq \left(|\psi(z)|+\frac{1-|\psi(z)|^2}{1-r^2}.\frac{r}{\cos r}\right)|g'(z)|\quad(|z|=r<1).
\end{equation}
Letting $|\psi(z)|=\beta$ $(0\leq \beta\leq1)$, the last inequality \eqref{abs f prime} becomes
\begin{equation*}
    |f'(z)|\leq \left(\beta+\frac{1-\beta^2}{1-r^2}.\frac{r}{\cos r}\right)|g'(z)|\quad(|z|=r<1).
\end{equation*}
Now we define
\begin{equation*}
  h(r,\beta):=\beta+\frac{1-\beta^2}{1-r^2}.\frac{r}{\cos r}\quad(0\leq \beta\leq1, 0<r<1).
\end{equation*}
To obtain $r_1$ it is enough to consider it as follows
\begin{equation*}
  r_1=\max\{r\in[0,1):h(r,\beta)\leq1, \forall \beta\in[0,1]\}.
\end{equation*}
Thus
\begin{align*}
  h(r,\beta)\leq1&\Leftrightarrow \frac{\beta(1-r^2)\cos r+(1-\beta^2)r}{(1-r^2)\cos r}\leq1 \\
  &\Leftrightarrow \beta(1-r^2)\cos r+(1-\beta^2)r\leq (1-r^2)\cos r\\
  &\Leftrightarrow 0\leq (1-\beta)(1-r^2)\cos r-(1-\beta^2)r\\
  &\Leftrightarrow 0\leq (1-r^2)\cos r-(1+\beta)r=:k(r,\beta).
\end{align*}
Since $\frac{\partial}{\partial\beta}k(r,\beta)=-r<0$, therefore we see that the above function $k(r,\beta)$ gets its minimum value in $\beta=1$, namely
\begin{equation*}
  \min\{k(r,\beta):\beta\in[0,1]\}=k(r,1)=:k(r),
\end{equation*}
where
\begin{equation*}
  k(r)= (1-r^2)\cos r-2r\quad(0<r<1).
\end{equation*}
Because $k(0)=1>0$ and $k(1)=-2<0$, therefore we conclude that there exists a $r_1$ such that $k(r)\geq0$ for all $r\in[0,r_1]$ where $r_1$ is the smallest positive
root of the Eq. \eqref{eq 1} and concluding the proof.
\end{proof}
Since the identity function $g(z)=z$ belongs to the class $\mathbf{S}^*_c$, thus letting $g(z)=z$ in the above Theorem \ref{th1} we get the following result. Indeed, we improve the bound $\sqrt{2}-1$ in the Theorem B.
\begin{corollary}
  If the function $f\in\mathcal{H}$ satisfies the conditions
  \begin{equation*}
    f(0)=0=f'(0)-1\quad{and}\quad |f(z)|<1,
  \end{equation*}
  then 
  \begin{equation*}
    |f'(z)|\leq 1\quad (|z|\leq 0.391389).
  \end{equation*}
\end{corollary}


\end{document}